\documentclass[12pt]{article}
\usepackage{amsmath}
\usepackage{amssymb}

\newtheorem{thm}{Theorem}[section]
\newtheorem{lem}{Lemma}[section]
\newtheorem{df}{Definition}[section]

\numberwithin{equation}{section}

\def\bbR{{\mathbb R}}
\def\bbZ{{\mathbb Z}}
\def\bbN{{\mathbb N}}
\def\bbQ{{\mathbb Q}}
\def\bbE{{\mathbf E}}
\def\bbP{{\mathbf P}}

\newcommand{\dd}{{\lfloor d/2 \rfloor}}

\def\eps{\varepsilon}
\let\phi=\varphi
\def\qed{\hfill $\square$}

\begin{document}

\title{The shape theorem for the frog
model with random initial configuration}

\author{O.S.M.~Alves$^{~1}$ \and
 F.P.~Machado$^{~2}$ \and
 S.Yu.~Popov$^{~2,3}$ \and
 K.~Ravishankar$^{~4}$ }

\maketitle

{\footnotesize
\noindent $^{~1}$Instituto de Matem\'atica e Estat\'\i stica,
Universidade Federal de Goias, Campus Samambaia, Caixa Postal 131,
CEP 74001--970, Goi\^ania GO, Brasil

\noindent e-mail: \texttt{oswaldo@mat.ufg.br}

\noindent $^{~2}$Instituto de Matem\'atica e Estat\'\i stica,
Universidade de S\~ao Paulo, Rua do Mat\~ao 1010, CEP 05508--900,
S\~ao Paulo SP, Brasil

\noindent e-mails: \texttt{fmachado@ime.usp.br, popov@ime.usp.br}

\noindent $^{~3}$Institute for Problems of Information Transmission
of Russian Academy of Sciences, Moscow, Russia

\noindent $^{~4}$Department of Mathematics, SUNY-New Paltz, NY
12561, USA

\noindent e-mail: \texttt{ravi@mcs.newpaltz.edu}

}

\begin{abstract}
We prove a shape theorem for a growing set of simple
random walks on $\bbZ^d$, known as frog model. The dynamics of
this process is described as follows: There are active particles,
which perform independent discrete time SRWs, and sleeping
particles, which do not move. When a sleeping particle is hit by
an active particle, the former becomes active as well.
Initially, a random number of particles is placed into each site.
At time~$0$ all particles are sleeping, except for those
placed at the origin. We prove that the set of all sites
visited by active particles, rescaled by the elapsed time,
converges to a compact convex set.
\\[.3cm] {\bf Keywords:} frog model, shape theorem, simple random walk
\end{abstract}

\section{Introduction}
In this note we study a discrete time
particle system in $ \bbZ^d $ named frog
model. In this model there are active particles, which move as
independent simple random walks (SRWs) on~$\bbZ^d$, and
sleeping particles, which do not move until activated.
At time zero there is a random number~$\eta(x)$
of particles at each site~$x$ of the
lattice, where $\{\eta(x), x\in\bbZ^d\}$ are i.i.d.,
and all the particles are sleeping except for those
that might be placed at the origin. Those active particles
start to perform a discrete time SRW.
  From then on when an active particle jumps on a
sleeping particle, the latter wakes up and starts jumping
independently, also performing a SRW.
If the origin is initially occupied then
the number of active
particles grows to infinity as active particles jump on sites that
have not been visited before, awakening the particles that are
sitting there. Let us underline that the active particles do not
interact with each other and there is no
``one-particle-per-site'' rule.

The frog model can be viewed as a model for describing information
spreading. The original idea is that every active particle has
some information and it shares that information with a sleeping
particle at the time the former jumps on the latter.
Particles that have the information move freely helping in the
process of spreading information. The model that we deal with in
this paper is a discrete-time version of that proposed by
R.~Durrett (1996, private communication),
who also suggested the term ``frog model''.

The first published result on
this model is due to Telcs, Wormald~\cite{TW}, where it was
referred to as the ``egg model''. They proved that, starting
from the one-particle-per-site initial configuration, the origin
will be visited infinitely often a.s. Popov~\cite{P}
proved that the last result holds in dimension $d\geq 3$
for the initial configuration with
a sleeping particle (or ``egg'') at each $x\neq 0$ with probability
$\alpha/\|x\|^2$, $\alpha$ being a large positive constant.
In Alves et al.~\cite{OFS2} a modification of the present
model was studied from the point of view of extinction and
survival. The difference of the model of~\cite{OFS2} from the
model of this paper is that in the former active particles may
disappear on each step. Recently A.~Ramirez and
V.~Sidoravicius communicated to us that they are working
on a continuous-time analog of this model, and that they
have proved some results such as shape theorem and
convergence to the product of Poissons.

In Alves et al.~\cite{OFS1}
it was proved that, starting from the one-particle-per-site
 initial configuration, the set of the
original positions of all active particles, rescaled by the
elapsed time, converges to a nonempty compact convex set.
In the present paper we generalize the main result of~\cite{OFS1}
to the case of random initial configuration. It turns out that
for the case when the initial configuration
contains empty sites, this generalization is nontrivial.

Now we define the model in a formal way.
Let $(\bbN^{\bbZ^d},{\cal B}_1,\nu)$
be a probability space,
where ${\cal B}_1$ is the product sigma
algebra and $\nu$ is the translation invariant
product measure determined by the distribution of
$ \{ \eta(x): x \in \bbZ^d \} $.

For each $\omega\in \bbN^{\bbZ^d}$ let $\{ {(S_{n,k}^x
(\omega))}_{n\in \bbN} , x\in \bbZ^d, 1\le k \le \omega(x) \}$
be the independent simple random walks which are executed by the
particles in $\omega$ when they are activated. We define
$S_{0,k}^x (\omega) =x$, for all $x \in \bbZ^d \hbox{ and } 1 \le k
\le \omega(x)$.  Denote by $\Omega^{\omega}$ the path
space of the trajectories of the random walks starting from the
initial configuration $\omega$ and by $\bbP_{\omega}$ the
corresponding path space measure.  Let $\bbP$ be the measure on
$\Omega = \prod_{\omega \in \bbN^{\bbZ^d}} (\omega \times
\Omega^{\omega} )$ obtained by taking the base measure on
$\bbN^{\bbZ^d}$ to be the product measure~$\nu$
and the conditional measure $\bbP [~\cdot
\mid \omega] = \bbP_{\omega}$. For each $ \omega \in \bbN^{\bbZ^d}$,
let
\[
t(x,z)(\omega) = \min \{ n : S^x_{n,k} (\omega) =z
\hbox{ for some } k, 1\le k \le \omega(x)\}
\]
and
\begin{equation}
\label{Tht}
T(x,z)(\omega) = \inf \Big\{\sum_{i=1}^m t(x_{i-1}, x_i)(\omega)\Big\},
\end{equation}
where the infimum is taken over all the finite sequences
$x=x_0,x_1,\ldots,x_m=z$.
Note that $t(x,z)(\omega) = T (x,z)(\omega) = \infty$
when $\omega(x)=0$. Note also that for $d\geq 3$,
$t(x,z)(\omega)=\infty$ with positive probability even
when $\omega(x)\geq 1$.

Let us define the set of sites which were visited by active
 particles up to time~$n$,
provided that initially the active particles were only in~$ x $.
Namely,
\[
 \xi_n^x (\omega) = \{ y \in \bbZ^d : T(x,y)(\omega) \le n \}.
\]
We are mostly concerned with $\xi_n:=\xi_n^0$
and its asymptotic behaviour.
In order to analyze that behaviour, define
\[
 {\bar \xi}^x_n (\omega) = \{ y + ({-1/2}, {1/2}]^d:
           y \in \xi^x_n(\omega)\}\subset\bbR^d,
\]
and ${\bar\xi}_n := {\bar\xi}^0_n$.

The main result of this paper is the following
\begin{thm}
\label{maintheo}
For any dimension $d\geq 1$
there is a nonempty convex set
$ {\mathcal A} = {\mathcal A}(d,\nu) \subset
\bbR^d $ such that for $\nu$-almost all initial
configurations~$\omega$, conditioned
on $\{\eta(0)\geq 1\}$, we have for any $0 < \eps < 1$
\[
(1-\eps){\mathcal A} \subset
\frac{{\bar \xi}_n}{n} \subset (1+\eps){\mathcal A}
\]
for all~$n$ large enough $\bbP_{\omega}$-a.s.
\end{thm}

Note that, although Theorem~\ref{maintheo}
establishes the existence of the
asymptotic shape~${\mathcal A}$, it is difficult to
identify exactly this shape. Of course, ${\mathcal A}$ is
symmetric and ${\mathcal A}\subset {\mathfrak D}$, where,
denoting $\|x\|_1=|x^{(1)}|+\cdots+|x^{(d)}|$,
\[
{\mathfrak D} = \{x\in \bbR^d :
     \|x\|_1 \leq 1\}.
\]
Also, note that if the initial configuration is augmented (i.e.\
some new particles are added), then the asymptotic shape (when it
exists) augments as well. In the paper~\cite{OFS1}
it was shown that if the
initial configuration is constructed by adding~$m$ particles to each
site and~$m$ is large enough, then the limiting
shape~${\mathcal A}$ contains some pieces of the
boundary of~${\mathfrak D}$ (a ``flat edge'' result). Now we show
that if the distribution of~$\eta$ is heavy-tailed enough, then
the limiting shape~${\mathcal A}$ coincides with~${\mathfrak D}$ (a
``full diamond'' result).

\begin{thm}
\label{full_diamond}
Suppose that for some positive $\delta<d$ and
for all~$n$ large enough we have
\begin{equation}
\label{fd}
 \bbP[\eta(x) \geq n] \geq (\log n)^{-\delta}.
\end{equation}
Then, Theorem~\ref{maintheo} is verified
with ${\mathcal A} = {\mathfrak D}$.
\end{thm}

\section{Proofs}
\noindent
{\it Proof of Theorem~\ref{maintheo}.}

\noindent
\underline{Step 1.} First, we state a few
preparatory results, which concern mainly the tails of
random variables $T(\cdot,\cdot)$.
Let us begin by recalling a technical fact from~\cite{OFS1}.
\begin{lem}
\label{expdecay}
Suppose that $\omega(x)=1$ for all~$x\in\bbZ^d$.
For all $ d \ge 1 $ and $ x_0 \in \bbZ^d $ there exist positive
finite constants $ \alpha_1=\alpha_1(x_0,d) $ and
$ \beta_1 = \beta_1(d)$ such that
\[
 \bbP [ T(0,x_0) \ge m ] \le \alpha_1 \exp\{- m^{\beta_1}\}
\]
for all~$m$.
\end{lem}

\noindent
{\it Proof.} Here we give only the main ideas of the proof,
as it was given in full detail in~\cite{OFS1} (Theorem~3.2).

\smallskip
\noindent
{\bf 1.} The case $d\geq 4$.
Pick $ n \ge \|x_0\|^2 $, where $\|\cdot\|$ is the
Euclidean norm, and fix some~$\eps$ such that $0<\eps<\frac{1}{2(d-2)}$.
Define for $1\leq i \leq \dd$ the sets
$\mathcal{D}^{(n)}_{i, \eps} := \{x\in\bbZ^d :
                     \|x\| \le i n^{{1/2}+ \eps}\}$.

\noindent
{\bf 1.1.} Consider the trajectory of the initial particle until
time~$n$. With overwhelming probability it
stays in $\mathcal{D}^{(n)}_{1, \eps}$
all that time and awakens at least $O(n^{1-\eps})$ sleeping particles.
To see why the last claim is true,
divide the time interval $ [0,n] $ into
 $ n^{\eps} $ disjoint subintervals of size $n^{1-\eps}$.
During a fixed subinterval, the expected size of the corresponding
subrange is of order $n^{1-\eps}$ (this can be seen by considering all the
sites~$y$ which are at distance
at most $\sqrt{n^{1-\eps}}$ from the position
of the initial particle at the beginning of the time subinterval,
estimating the probability that~$y$ is hit and summing over~$y$;
or just by applying directly the known results about the
size of the range of
SRW, see e.g.\ Hughes~\cite{Hughes}, pp.~333, 338). As the size of the
subrange is certainly not greater than $n^{1-\eps}$, with
probability bounded away from~$0$ it will be $O(n^{1-\eps})$
(this follows from the following fact:
For any random variable~$X$ with $0\leq X\leq a$ a.s.\ and
$\bbE X\geq b$ it is true that $\bbP[X\geq b/2]\geq b/(2a)$).
So, with overwhelming probability the size of at least one of
the subranges will be of order $n^{1-\eps}$.

\noindent
{\bf 1.2.} Pick $n^{1-\eps}$ particles in $\mathcal{D}^{(n)}_{1, \eps}$
 which were awakened by the initial particle by time~$n$.
Divide them into $n^{\eps}$ disjoint groups of size $n^{1-2\eps}$.
Observe the particles of a fixed group during $n^{1+2\eps}$ time
units after their activation. Let~$\zeta$ be the number of sites
in $\mathcal{D}^{(n)}_{2, \eps}\setminus \mathcal{D}^{(n)}_{1, \eps}$
visited by the particles of that group along the time interval mentioned
above.
Clearly, $\zeta\leq n^{1-2\eps}\times n^{1+2\eps} = n^2$
and the direct computation (for each
$y\in \mathcal{D}^{(n)}_{2, \eps}\setminus
          \mathcal{D}^{(n)}_{1,\eps}$, using that the $n^{1-2\eps}$
particles from the group are independent,
we compute a lower bound on the probability that at least
one particle of the group hits~$y$, and then sum over~$y$) shows that
$\bbE\zeta = O(n^2)$, so with probability bounded away from~$0$,
$\zeta$ is of order~$n^2$. Considering now all the~$n^\eps$ groups
and using the independence again,
we obtain that with overwhelming probability there will be $O(n^2)$
particles in
$\mathcal{D}^{(n)}_{2, \eps}\setminus \mathcal{D}^{(n)}_{1, \eps}$
by time $n+n^{1+2\eps}$.

\noindent
{\bf 1.3.} Now, proceeding in the same spirit, we use those~$O(n^2)$
particles from
$\mathcal{D}^{(n)}_{2, \eps}\setminus \mathcal{D}^{(n)}_{1, \eps}$
 to awaken~$O(n^3)$ particles
in $\mathcal{D}^{(n)}_{3, \eps}\setminus \mathcal{D}^{(n)}_{2, \eps}$,
and so on, to get finally $O(n^{\dd})$
active particles in
$\mathcal{D}^{(n)}_{\dd, \eps}\setminus \mathcal{D}^{(n)}_{\dd-1, \eps}$
 at time $n+O(n^{1+2\eps})$ with overwhelming probability.

\noindent
{\bf 1.4.} Considering those~$O(n^{\dd})$
 particles in
 $\mathcal{D}^{(n)}_{\dd, \eps}\setminus \mathcal{D}^{(n)}_{\dd-1, \eps}$
and waiting $n^{1+2\eps}$ units of time more, one gets that with
overwhelming probability at least one of those particles will hit~$x_0$
by time $m_d:=n+O(n^{1+2\eps})$.

\smallskip
\noindent
{\bf 2.} The case $d=3$. Again, let $n\geq \|x_0\|^2$.
Take $\eps<1/4$ and consider the $n^{1-\eps}$
particles in $\mathcal{D}^{(n)}_{1, \eps}$
awakened by the initial particle until
time~$n$. Now, as in the last step of the argument for $d\geq 4$,
until time $m_3:=n+O(n^{1+2\eps})$ with overwhelming probability
at least one of those particles will hit~$x_0$.

\smallskip
\noindent
{\bf 3.} The case $d=2$. This case
 is treated analogously to the case $d=3$.
That is, first,
dividing the time interval $[0,n]$
into $n^{1/2}$ subintervals of size
$n^{1/2}$ (i.e., taking $\eps=1/2$) and using the
fact that the expected size of the range of SRW grows like
$O(n/\log n)$, we
 prove that with large probability the original
particle will awaken $O(n^{1/2}/\log n)$ sleeping particles
in the ball of radius~$n$ until the moment~$n$,
where $n\geq\|x_0\|^2$.
Considering the independent random walks of
those particles until time $m_2:=n+O(n^2)$,
 we get the result.

\smallskip
\noindent
{\bf 4.} The case $d=1$. This was not treated in Theorem~3.2
of~\cite{OFS1}, but anyway it is quite analogous to the cases $d=2,3$.
First, by time~$n$ we will have $n^{1/4}$
active particles situated not farther
than~$n$ from the origin. Then, waiting until $m_1:=n+n^2$ one gets the
result.  \qed

\medskip

For $d\geq 3$ denote $\eps_d=(6(d-2))^{-1}$; clearly, when doing the proof
of Lemma~\ref{expdecay} in dimension~$d\geq 3$,
one may fix $\eps=\eps_d$.
  From the above proof one can deduce that there exist
deterministic constants~$h_d$ (which depend only on dimension)
such that $m_d\leq h_d n^{1+2\eps_d}$, $d\geq 3$, $m_d\leq h_dn^2$,
$d=1,2$. This means that, in order to obtain an upper bound on
$\bbP[T(0,x_0)\geq m]$, one should
follow the steps of the proof of Lemma~\ref{expdecay}
with $n=(m/h_d)^{1/(1+2\eps)}$, $d\geq 3$, or $n=(m/h_d)^{1/2}$,
$d=1,2$. Keeping this in mind and denoting
\[
{\cal R}_{\omega} (B) = \frac{1}{|B|} \sum_{x\in B}
                  {\bf 1}_{\{\omega(x)\geq 1\}},
\]
where $\omega$ is the initial configuration and~$B$ is a finite subset
of~$\bbZ^d$, consider the following
\begin{df}
\label{goodeta}
Let $p_1=\bbP[\eta(0)\geq 1]$, and $m$ be any positive integer.
A fixed initial configuration~$\omega$ is called $m$-good if
\begin{itemize}
\item $d\geq 4$ and
  \begin{itemize}
   \item for any ball~$B$ of radius $n_d(m)^{(1-\eps_d)/2}$ which
   is fully inside $\mathcal{D}^{(n_d(m))}_{1, \eps_d}$ we have
   ${\cal R}_{\omega}(B)\geq p_1/2$;
   \item we have
   ${\cal R}_{\omega}(\mathcal{D}^{(n_d(m))}_{i, \eps_d}\setminus
      \mathcal{D}^{(n_d(m))}_{i-1, \eps_d})\geq p_1/2$
       for all $i=2,\ldots,\dd$;
  \end{itemize}
\item $d=3$, and
 for any ball~$B$ of radius $n_3(m)^{(1-\eps_3)/2}$ which
   is fully inside $\mathcal{D}^{(n_3(m))}_{1, \eps_3}$ we have
   ${\cal R}_{\omega}(B)\geq p_1/2$;
\item $d=1,2$, and for any ball~$B$ of radius $n_d(m)^{1/4}$
situated not farther than~$n_d(m)$ from the origin, we have
${\cal R}_{\omega} (B) \geq p_1/2$,
\end{itemize}
where $n_d(m)=(m/h_d)^{1/(1+2\eps_d)}$,
$d\geq 3$, or $n_d(m)=(m/h_d)^{1/2}$,
$d=1,2$, and the notation $\mathcal{D}^{(n)}_{i, \eps_d}$ is from
the proof of Lemma~\ref{expdecay}.
\end{df}

\begin{lem}
\label{expdecay2}
For all $ d \ge 1 $ and $ x_0 \in \bbZ^d $ there exist positive
finite constants $ \alpha_2=\alpha_2(d,p_1) $ and
$ \beta_2 = \beta_2(d,p_1)$ such that if~$\omega$ is $m$-good
and $\|x_0\|^2\leq n_d(m)$ (the notation $n_d(m)$ is
 from Definition~\ref{goodeta}), then
\[
 \bbP_{\omega} [ T(0,x_0) \ge m ] \le \alpha_2 \exp\{- m^{\beta_2}\}.
\]
\end{lem}

\noindent
{\it Proof.} From Definition~\ref{goodeta} it follows that,
by following the steps of the proof of Lemma~\ref{expdecay}
one can prove Lemma~\ref{expdecay2}. For example, in the part 1.1
the expected amount of sleeping
particles activated during a fixed time
subinterval differs from the expected size of the subrange only
by a constant factor (depending only on~$p_1$ and~$d$) when~$\omega$
is $m$-good. In the part 1.2, $\bbE\zeta$ will be of the same order
for $m$-good $\omega$ and for one-particle-per-site $\omega$.
In general, each time that one is computing the expected number
of newly awakened particles in the proof of Lemma~\ref{expdecay},
considering $m$-good initial configuration instead of
one-particle-per-site initial configuration costs only a
constant factor.
Another observation, which is crucial for the subsequent
discussion, comes after examining the proof of
Lemma~\ref{expdecay}: As long as $\|x_0\|^2\leq n$, the estimates
that one gets on $\bbP[T(0,x_0)>m_d]$ are uniform in~$x_0$.
\qed

\medskip

Define $\hat n_d(m)$ to be equal to $n_d(m)^{(1-\eps(d))/3}$,
$d\geq 3$, and to $n_d(m)^{1/8}$, $d=1,2$.

\begin{lem}
\label{expdecay3}
For all $ d \ge 1 $ there exist positive
finite constants $ \alpha_3=\alpha_3(d,p_1) $ and
$ \beta_3 = \beta_3(d,p_1)$ such that
\[
 \bbP [\eta \mbox{ is $m$-good} \mid \eta(x_1)=\cdots=\eta(x_{k_0})=0]
    \geq 1 - \alpha_3 \exp\{- m^{\beta_3}\}
\]
for any fixed collection of sites $x_1,\ldots,x_{k_0}\in \bbZ^d$,
$k_0\leq \hat n_d(m)^d$.
\end{lem}

\noindent
{\it Proof.} Note that $\hat n_d(m)^d$ is small relative to the
sizes of sets which were considered in Definition~\ref{goodeta}.
As $\eta(x)$ are i.i.d., it is straightforward
to prove this fact by using the large deviations technique.
\qed

\medskip

In the sequel we will need the following basic fact which
is stated without proof here.
\begin{lem}
\label{elem_fact}
Let $(X_i, i\geq 1)$ be a sequence of nonnegative random
variables, and there exist $\phi_n$, $n=1,2,\ldots$, such that
$\sup_i \bbP[X_i\geq n]\leq \phi_n$ and $\sum_n \phi_n < \infty$.
Then $X_n/n \to 0$ a.s.
\end{lem}

\medskip
\noindent
\underline{Step 2.} Now the goal is to verify the
conditions of Kingman-Liggett subadditive ergodic theorem \cite{K,L}
for the sequence of random variables $Y(m,n)$ defined below.

For a fixed $x\in \bbZ^d$ and an $\omega \in \bbN^{\bbZ^d}$
satisfying the condition  $\omega (0) \geq 1$
define a sequence of positive integers
${\{ v^x_k\}}_{k=0}^\infty$ as follows:
\begin{eqnarray*}
v^x_0 &=& 0,\\
v^x_{k+1} & = & \min \{n > v^x_k : \omega(n x) \geq 1 \}.
\end{eqnarray*}

In words, $v^x_k x$, $k=0,1,2,\ldots$, are the sites on the ray
$\{n x, n\geq 0\}$ which are
occupied in the initial configuration.
Clearly, for any $x\in\bbZ^d$ and $i=1,2,\ldots$ it holds that
$\bbP[v^x_i-v^x_{i-1}=k]=p_1(1-p_1)^{k-1}$ ($p_1$ is from
Definition~\ref{goodeta}).
 Now, for $m,n\geq 0$ consider the collection of random variables
\[
  Y(m,n) = T(v^x_mx,v^x_nx).
\]
It is important to observe that
the random variables $ \{T(x,y): x,y \in \bbZ^d\}$
are subadditive in the sense that for {\it all\/} $\omega$,
\begin{equation}
\label{dfg}
T(x,z) \le T(x,y) + T(y,z) \hbox{ for all } x,y,z \in \bbZ^d.
\end{equation}
Here is the explanation.
Fix the initial configuration~$\omega$.
If the site~$y$ was empty in the initial configuration, then
$T(y,z)=\infty$ and~(\ref{dfg}) follows.
Now, suppose that $\omega(y)\geq 1$.
If site~$z$ is reached before site~$y$,
then~(\ref{dfg}) is evident. If that does not happen, recall that
the random variables $T(y,z)$, $y,z\in\bbZ^d$ are constructed
using the same collection of the random variables $S^x_{n,k}$,
i.e., each particle follows the same trajectory as soon as it wakes up.
So the process departuring from only site~$y$
awakened (the one which gives the
passage time $T(y,z)$) is coupled with the original
process (i.e., that started from~$x$), and for the latter one may have
other particles awakened at time $T(x,y)$ besides that from~$y$.
Consequently, from~(\ref{Tht}) it follows that $T(x,z)-T(x,y)$,
which is the remaining time to reach site~$z$ for the original
process, is less than or equal to $T(y,z)$, thus proving~(\ref{dfg}).
The equation~(\ref{dfg}) shows that for {\it all\/}
initial configurations~$\omega$,
\begin{equation}
\label{dfg1}
Y(m,k) \le Y(m,n) + Y(n,k) \hbox{ for all } m,n,k \in \bbN.
\end{equation}

Now, let us verify the conditions of
Kingman-Liggett subadditive ergodic theorem for
the random variables $Y(m,n)$.

Since $\nu$ is a product measure we can
 easily establish the following
result:
Given nonnegative integer numbers $m_1, m_2, \ldots ,m_n, p$,
 the joint distributions of random variables
$\{ Y(m_1,m_2), \ldots ,Y(m_{n-1},m_{n})\}$
 and
$\{ Y(m_1+p ,m_2+p ), \ldots,Y(m_{n-1}+p ,m_{n}+p )\} $
  are equal. From this we obtain the stationarity
  of the sequence ${\{ Y((n-1)k,n k) \}}_{n\in \bbN}$
for each $k \in \bbN$.
Ergodicity of the sequence follows from the observation that
the events $\{T(v^x_{n_1 k}x,v^x_{(n_1 +1)k}x) \le a \}$  and
$\{ T(v^x_{n_2 k}x, v^x_{(n_2 + 1) k}x ) \le b \}$
are independent provided
$a+b \le \| (n_1 - n_2)k x \|_1$.

Next we prove that  $\bbE T(0,v^x_1 x) < \infty$.
Using Lemmas~\ref{expdecay2} and~\ref{expdecay3},
for~$y$ such that $\|y\|\leq \hat n_d(m)^{1/2}$
(note that $\hat n_d(m) < n_d(m)$) and for any collection
of sites $z_1,\ldots,z_{k_0}\in\bbZ^d$,
$k_0\leq \hat n_d(m)^d$, we have,
denoting $\bbP^*[\,\cdot\,]=
 \bbP[~\cdot \mid \eta(z_1)=\cdots=\eta(z_{k_0})=0]$,
\begin{eqnarray}
 \bbP^*[T(0,y) \geq m]
&\leq&  \bbP^*[T(0,y)\geq m \mid \eta
        \mbox{ is $m$-good}] +
       \bbP^*[\eta \mbox{ is not $m$-good}]\nonumber\\
&\leq&  \alpha_2\exp\{-m^{\beta_2}\} +
            \alpha_3\exp\{-m^{\beta_3}\}.\label{matozh}
\end{eqnarray}
Using~(\ref{matozh}), we get
\begin{eqnarray}
  \bbP[T(0,v^x_1x)\geq m]  &=&
    \sum_{k=1}^\infty \bbP[v^x_1=k] \bbP[T(0,k x)\geq m\mid v^x_1=k]
                                               \nonumber\\
 &\leq& \sum_{k\leq \|x\|^{-1}\hat n_d(m)^{1/2}} \bbP[v^x_1=k]
           \bbP[T(0,k x)\geq m\mid v^x_1=k]  \nonumber\\
 &&{}+ \sum_{k > \|x\|^{-1}\hat n_d(m)^{1/2}} p_1(1-p_1)^{k-1}\nonumber\\
  &\leq & \alpha_2\exp\{-m^{\beta_2}\} +
            \alpha_3\exp\{-m^{\beta_3}\} \nonumber\\
           &&{} +(1-p_1)^{\|x\|^{-1}\hat n_d(m)^{1/2}}, \label{Evk}
\end{eqnarray}
so $\bbE Y(0,1) = \sum_{m\geq 1} \bbP[T(0,v^x_1x)\geq m] < \infty$.

Thus, we have verified the conditions of Kingman-Liggett
subadditive ergodic theorem for the sequence
$\{ Y(m,n) = T(v^x_m x, v^x_n x)\}$.
Therefore, one gets that there exists $\mu'(x)$ such that
\begin{equation}
\label{mu'}
\lim_{n\to \infty} \frac{ Y(0,n)}{n} = \mu'(x) <\infty \qquad
\mbox{$\bbP$-a.s.}
\end{equation}

\medskip
\noindent
\underline{Step 3.} Now the goal is to pass from the random
sequence $(v^x_k x, k\geq 0)$ to the whole ray $(k x, k\geq 0)$.

Denote $\mu(x):=p_1\mu'(x)$. Let $n \in \bbN$ satisfy
$v^x_{k(n)-1} \leq n < v^x_{k(n)}$. Then, by the
subadditivity,
\begin{equation}
\label{u***}
\frac{T(0,n x)}{n} \le \frac{T(0,v^x_{k(n)}x)}{k(n)}\frac{k(n)}{n} +
\frac{T(v^x_{k(n)}x, n x)}{n}.
\end{equation}
Note that, since $\nu$ is a product measure,
$v^x_{k(n)}- n$ has the same distribution
as $v^x_1$.
This means that the upper bound on $\bbP[T(v^x_{k(n)}x, n x)\geq m]$
is also given by~(\ref{Evk}), so by Lemma~\ref{elem_fact}
  we have that
$n^{-1}T(v^x_{k(n)}x,n x) \to 0$, $\bbP$-a.s.
Together with~(\ref{mu'}) and
the fact that $k(n)/n \to p_1$ $\bbP$-a.s.,
this shows that one gets from~(\ref{u***}) that
\begin{equation}
\label{lsup}
 \limsup_{n\to\infty} \frac{T(0, n x)}{n} \le \mu(x)
   \qquad \mbox{$\bbP$-a.s.}
\end{equation}

Now, let us prove that
\begin{equation}
\label{linf}
 \liminf_{n\to\infty} \frac{T(0, n x)}{n} \ge \mu(x)
   \qquad \mbox{$\bbP$-a.s.}
\end{equation}
For a fixed site $y\in\bbZ^d$, $y\neq 0$,
let $U_y$ be a random variable defined in the following way.
At the moment $T(0,y)$ consider the active particle situated
in~$y$ (if there are several such particles, choose one of them by
randomization). Let $Z^y_k, k\geq T(0,y)$ be the subsequent walk
of this particle. We define
$\tau = \min\{k: \eta(Z^y_k)\geq 1\}$ and so $U_y:=\tau-T(0,y)$
is the time that the particle needs to travel from~$y$ to
some site that initially contained sleeping particles.
Clearly, $U_y=0$ corresponds to the case $\eta(y)\geq 1$.
Now, let $\hat\zeta_n = Z^{nx}_{T(0, n x)+U_{nx}}$;
note that $\eta(\hat\zeta_n)\geq 1$.
By subadditivity one can write
\begin{equation}
\label{u*}
 T(0,v^x_{k(n)}x) \leq T(0, n x ) + U_{n x}
                  + T(\hat\zeta_n, v^x_{k(n)}x).
\end{equation}
To proceed, we need to get some good bounds on the tails
of~$U_{n x}$ and $T(\hat\zeta_n, v^x_{k(n)}x)$.
Note that there exists a positive constant~$\beta$
which is not dependent on the dimension, such that
$\bbP[|\{Z^{n x}_j, 0\leq j\leq k^{1/2}\}|\geq k^{1/4}] \geq \beta$,
for all $k\geq 1$. This shows that
\[
 \bbP[|\{Z^{n x}_j, 0\leq j\leq k\}|\geq k^{1/4}] \geq
   1 - (1-\beta)^{k^{1/2}}.
\]
So, as any site in $\{Z^{n x}_j, 0\leq j\leq k\}$ initially contained
a sleeping particle with probability~$p_1$, it holds that
\begin{equation}
\label{EUk}
 \bbP[U_{n x} \geq k] \leq (1-\beta)^{k^{1/2}} +
 (1-p_1)^{k^{1/4}}.
\end{equation}
Now, as $\|n x - \hat\zeta_n\| \leq U_n$, (\ref{EUk}) implies
that
\begin{eqnarray*}
\bbP[\|v^x_{k(n)}x - \hat\zeta_n\|\geq k] &\leq &
 \bbP[\|v^x_{k(n)}x - n x\|\geq k/2] +
  \bbP[\|n x - \hat\zeta_n\|\geq k/2]\\
  &\leq& (1-p_1)^{k/2} +
    (1-\beta)^{(k/2)^{1/2}} + (1-p_1)^{(k/2)^{1/4}}\\
    &=:&\psi(k),
\end{eqnarray*}
and, proceeding similarly to~(\ref{Evk}), we get that
\begin{equation}
\label{fff}
\bbP[T(\hat\zeta_n, v^x_{k(n)}x)\geq m] <
   \alpha_2\exp\{-m^{\beta_2}\} +
            \alpha_3\exp\{-m^{\beta_3}\} +
                    \psi(\|x\|^{-1}\hat n_d(m)^{1/2}).
\end{equation}

Thus, by~(\ref{EUk}) and~(\ref{fff}), dividing~(\ref{u*})
by~$n$ and applying~(\ref{mu'}) and
Lemma~\ref{elem_fact}, we get~(\ref{linf}).
Now, from~(\ref{lsup}) and~(\ref{linf}) we finally obtain that,
for all $x\in\bbZ^d$, there exists $\mu(x) \geq 0$ such that
\begin{equation}
\label{limtint}
 \lim_{n\to\infty}\frac{T(0,n x)}{n} \to \mu(x) \qquad
    \mbox{$\bbP$-a.s.}
\end{equation}

Notice that the equation~(\ref{limtint}) is already enough to
get the proof of Theorem~\ref{maintheo} in dimension~1 with
${\mathcal A}=[-(\mu(1))^{-1}, (\mu(1))^{-1}]$. So, from now on
we shall concentrate on the case $d\geq 2$.

\medskip
\noindent
\underline{Step 4.}
To proceed with the proof of Theorem~\ref{maintheo},
one has to assure that~$\xi_n$ grows
{\it at least\/} linearly.
For $y\in\bbR^d$ and $a>0$ denote $D(y,a) = \{x\in\bbR^d :
   \|x-y\|_1 \leq a\}$.
\begin{lem}
\label{smallball}
For all $d\geq 2$ there exist
constants $ 0 < \delta < 1 $, $\alpha_4>0$, $\beta_4>0$,
which depend only on the dimension, such that
\[
 \bbP[ D(x, {n \delta})
\subset {\bar \xi}_{n+T(0,x)} ] \geq
                     1-\alpha_4\exp\{-n^{\beta_4}\}
\]
for all~$n$ and~$x$, conditioned on $\{\eta(0)\geq 1\}$.
\end{lem}

\noindent
{\it Proof.} The proof of this fact follows the spirit of
the proof of Lemmas~4.2 and~4.3 of~\cite{OFS1}.
The difficulty that arises here is that, for~$\omega$
such that~$\omega(x)=0$, $T(x,y)(\omega)=\infty$ for any~$y$,
so a direct application of the method of~\cite{OFS1} does not work.
To get around that difficulty, we use the construction similar to
that used in Step~3. Recall the definition of random variable~$U_y$
 from Step~3, and define for $y\neq 0$,
$T^*(y,z) = U_y+T(\hat\zeta^y,z)$, where, similarly to Step~3,
$\hat\zeta^y$ is the site with $\eta(\hat\zeta^y)\geq 1$
at which the active particle which was
in~$y$ at time $T(0,y)$ arrived after~$U_y$ steps.
Now, let~$y$ be such that $\|x-y\|_1=n$, and let
$x=y_0,y_1,\ldots,y_n=y$ be a path connecting~$x$ to~$y$
such that for all~$i$, $\|y_i-y_{i-1}\|_1=1$;
note that $\|x-y_k\|_1=k$, $k=0,\ldots,n$.
Denote $Y_i = T^*(y_{i-1},y_i)$. From the above construction it
follows that
\[
T(0,y)-T(0,x) \leq \sum_{i=1}^n Y_i =
 \sum_{i=1}^{n^{1/2}}\sigma_i,
\]
where
\[
\sigma_i = \sum_{j:i+jn^{1/2}\leq n} Y_{i+jn^{1/2}},
\]
$i=1,\ldots,n^{1/2}$.
Analogously to~(\ref{EUk}) and~(\ref{fff}), it is not difficult
to get that for some constants $\alpha_5,\beta_5>0$ and for all~$m$
\begin{equation}
\label{u****}
\bbP[Y_i\geq m] \leq \alpha_5\exp(-m^{\beta_5}),
\end{equation}
uniformly in~$x,y$ and paths connecting them.
Consider the event $B=\{Y_i<n^{1/2}/2, i=1,\ldots,n\}$;
clearly, from~(\ref{u****}) it follows that for
some $\alpha_6,\beta_6>0$, $\bbP[B]\geq 1-\alpha_6\exp(-n^{\beta_6})$.
Note that, if $|i-j|\geq n^{1/2}$, then the random variables
$Y_i {\bf 1}_B$ and $Y_j {\bf 1}_B$ are independent, since if the
event~$B$ occurs, then the random variables~$Y_i$ and~$Y_j$
depend on disjoint sets of random walks. So, when~$B$ occurs,
each~$\sigma_i$ is a sum of independent random variables.
Although~$Y_i$ are not identically distributed and we cannot
guarantee the existence of moment generating function
for~$Y_i$, the condition~(\ref{u****}) allows us, by
applying Theorem~1.1 of Nagaev~\cite{Nagaev}, to obtain that
there exists~$\delta_0>0$ such that $\sigma_i \leq n^{1/2}/\delta_0$
with probability at least subexponentially high.
This shows that, with overwhelming probability,
$\sum_{i=1}^n Y_i \leq n/\delta_0$. Thus, if~$y$ is
at distance~$n$ from~$x$, with overwhelming probability
by time $T(0,x)+n/\delta_0$ it will be visited. Now,
if $d\geq 2$ and $0<\|x-y\|_1<n$, then there exists $z\in\bbZ^d$
such that $\|z-x\|_1=n$ and $\|z-y\|_1$ equals~$n$ or~$n+1$. Since
$T^*(x,y)\leq T^*(x,z)+T^*(z,y)$, the result follows with
$\delta=\delta_0/2$.
\qed

\medskip
\noindent
\underline{Step 5.} The next step is to prove that
$\mu(\cdot)$ can be extended to~$\bbR^d$ in such a way that
$\mu$ is a norm in~$\bbR^d$.
Let us extend the definition of $ T(x,y) $ to the
whole $ \bbR^d \times \bbR^d $ by defining
\[
 T(x,y) = \min\{n:y \in {\bar \xi}^{x_0}_n\},
\]
where $x_0\in \bbZ^d$ is such that $x\in (-1/2,1/2]^d + x_0$.
   From the fact $ T(0,n x) \ge n\|x\|_1 $
it follows that $ \mu(x) \ge \|x\|_1 $
for all $ x \in \bbZ^d$.

\begin{lem}
\label{normZd}
For any $a\in\bbR_+$, $x,y\in\bbZ^d$ we have
\begin{equation}
\label{norm1}
\mu(a x) := \lim_{n\to\infty} \frac{T(0, a n x)}{n} = a\mu(x)
    \qquad \mbox{$\bbP$-a.s.},
\end{equation}
and
\begin{equation}
\label{norm2}
\mu(x+y) \leq \mu(x)+\mu(y).
\end{equation}
\end{lem}

\noindent
{\it Proof.} First, the proof of~(\ref{norm1}) basically
repeats what was done on Step~3, so we omit it.
Let us turn to the proof of~(\ref{norm2}). Clearly,
instead of~(\ref{norm2}) it is sufficient to prove
that for any $x,y\in\bbZ^d$
\begin{equation}
\label{norm2'}
\mu'(x+y) \leq \mu'(x)+\mu'(y)
\end{equation}
with $\mu'(\cdot)$ defined by~(\ref{mu'}).

For fixed $n\geq 1$ define a random sequence $s_{k,n}$,
$k=0,1,2,\ldots,$ in the following way:
\begin{eqnarray*}
s_{0,n} &=& 0,\\
s_{k+1,n} &=& \min\{m>s_{k,n}:\eta(my+v^x_n x)\geq 1\},
\end{eqnarray*}
and let $z_{k,n}=s_{k,n}y+v^x_n x$. By the subadditivity,
\begin{equation}
\label{zzz}
T(0,v^{x+y}_n (x+y)) \leq T(0,v^x_n x) + T(v^x_n x, z_{n,n})
  + T(z_{n,n},v^{x+y}_n (x+y)).
\end{equation}
We have
\begin{eqnarray*}
v^{x+y}_n (x+y) - z_{n,n} &=& x\sum_{i=1}^n((v^{x+y}_i-v^{x+y}_{i-1})
   -(v^x_i-v^x_{i-1}))\\
  &&  {}+y\sum_{i=1}^n
    ((v^{x+y}_i-v^{x+y}_{i-1})-(s_{i,n}-s_{i-1,n})).
\end{eqnarray*}
Fix arbitrary $\eps>0$.
Now, $(v^{x+y}_i-v^{x+y}_{i-1})$, $(v^x_i-v^x_{i-1})$,
$(s_{i,n}-s_{i-1,n})$, $i=1,\ldots,n$, are i.i.d.\ random
variables, geometrically distributed with parameter $(1-p_1)$.
So, in the right-hand side of the above display we have
two sums of i.i.d.\ random variables satisfying the Cramer
condition and with zero mean. So, by the Large Deviation
Principle, with exponentially high probability
$\|v^{x+y}_n (x+y) - z_{n,n}\| \leq \eps(\|x\|+\|y\|)n$. By
Lemma~\ref{smallball}, with subexponentially high probability,
$T(z_{n,n},v^{x+y}_n (x+y)) \leq \eps\delta^{-1}(\|x\|+\|y\|)n$
(the site $z_{n,n}$ is random, so, to see that
Lemma~\ref{smallball} is applicable here, note that
$T(z_{n,n},v^{x+y}_n (x+y))$ has the same distribution
as $T(0,u)$, where~$u$ is a random site constructed as follows:
First, let~$u'$ be the $n$-th occupied site on the ray along~$(-y)$,
then, $u''$ is the $n$-th occupied site on the ray along~$(-x)$
beginning in~$u'$, and finally,
$u$ is the $n$-th occupied site on the ray along~$x+y$
beginning in~$u''$).
Note that $T(v^x_n x, z_{n,n})$ has the same distribution
as $T(0,v^y_n y)$, so $n^{-1}T(v^x_n x, z_{n,n}) \to \mu'(y)$
at least in probability. Dividing~(\ref{zzz}) by~$n$ and
taking the limit in probability, we get
\[
\mu'(x+y) \leq \mu'(x)+\mu'(y)+\eps\delta^{-1}(\|x\|+\|y\|);
\]
when $\eps\to 0$, we obtain~(\ref{norm2'}). \qed

\medskip

Now, let us show that $\mu(x)$ can be extended to the
whole~$\bbR^d$ in such a way that~(\ref{limtint}) still holds.
This is done in a standard way: First, if $x\in\bbQ^d$
(i.e., all the coordinates of~$x$ are rational), we define
$\mu(x):={\mu(m x)}/{m}$, where~$ m $ is the smallest
positive integer such that $ m x \in \bbZ^d $
(the fact that~(\ref{limtint}) still holds for $x\in\bbQ^d$
follows from~(\ref{norm1})).
Then, by using the fact that~$\mu$ is a norm in~$\bbQ^d$
(this follows from Lemma~\ref{normZd}), it is
extended to~$\bbR^d$. To prove that~(\ref{limtint}) holds
for any~$x\in\bbR^d$, one can proceed as follows:
By the subadditivity and Lemma~\ref{smallball},
if $y\in\bbQ^d$, then
$\limsup_{n\to\infty}n^{-1}T(0,n x)\leq \mu(y) +
\delta^{-1}\|x-y\|$, and
$\liminf_{n\to\infty}n^{-1}T(0,n x)\geq \mu(y) -
\delta^{-1}\|x-y\|$. Approximating~$x$ by vectors with
rational coordinates, one gets the result.

\medskip
\noindent
\underline{Step 6.}
Now everything is ready to finish the proof of Theorem~\ref{maintheo}.

Let $ { \mathcal A } := \{x \in \bbR^d : \mu(x) \le 1 \} $.
The following argument is standard (see e.g.~\cite{BG,D,DG}),
we keep it to preserve the self-containedness of the paper.

Denote $\eps'=(1-\eps)^{-1}-1$, and $\eps''=1-(1+\eps)^{-1}$.
To prove Theorem~\ref{maintheo}, it is enough to prove
that $ n{\mathcal A} \subset {\bar \xi}_{(1+\eps')n}$ and
${\bar \xi}_{(1-\eps'')n} \subset n{\mathcal A} $
for all~$ n $ large enough, $\bbP$-a.s.\
(and so $\bbP_\omega$-a.s.\ for $\nu$-almost all~$\omega$).

Choose a finite set
$ F:= \{x_1, \ldots, x_k\} \subset {\mathcal A} $ such that
$ \mu(x_i) < 1 $ for $i=1,\ldots,k$, and
(with~$\delta$ from Lemma~\ref{smallball})
$ {\mathcal A} \subset \cup_{i=1}^{k}D (x_i, \eps'\delta)$.
Notice that~(\ref{limtint}) implies that
$ n F \subset
{\bar \xi}_{n}$ for all $n$ large
enough, $\bbP$-a.s.
Now, Lemma~\ref{smallball} and Borel-Cantelli imply that $\bbP$-a.s.\
for all~$n$ large enough we have $D( n x_i, n\eps'\delta) \subset
{\bar \xi}^{ n  x_i}_{n \eps'}$, for all $i=1,2, \dots ,k$.
So  $ n{\mathcal A} \subset {\bar \xi}_{(1+\eps')n}$
 and this part of the proof is done.

Now, let us choose $ G := \{y_1, \ldots, y_k\} \subset
2{\mathcal A}\setminus {\mathcal A} $ in such a way that
$ 2{\mathcal A}\setminus {\mathcal A}
\subset \cup_{i=1}^k D(y_i, \eps''\delta )$.
Analogously, we get that
$ n G \cap {\bar \xi}_{n} =\emptyset$
 for all~$n$ large enough $\bbP$-a.s., and
that for all~$n$ large enough, if ${\bar \xi}_{(1-\eps'')n} \cap
(2{\mathcal A}\setminus {\mathcal A}) \neq \emptyset$, then
${\bar \xi}_{n} \cap n G \neq \emptyset$. This shows
that ${\bar \xi}_{(1-\eps'')n} \subset n{\mathcal A} $
for all~$ n $ large enough, $\bbP$-a.s., and so concludes
 the proof of Theorem~\ref{maintheo}. \qed

\medskip
\noindent
{\it Proof of Theorem~\ref{full_diamond}.} Denote
\[
{\mathfrak D}_n = \{x \in \bbZ^d : \|x\|_1 \leq n\}.
\]
Choose $\theta<1$ such that $\delta<\theta d$. Start the process
and wait until the moment~$n^\theta$. By Lemma~\ref{smallball}
there exist~$C_1$, $\gamma_0$, such that with probability at least
$1-\exp\{n^{\gamma_0}\}$ all the frogs
which were initially in the ball of
radius~$C_1n^\theta$ centered in~$0$ will be awake.
Clearly, the inequality~(\ref{fd}) implies that
\[
 \bbP[\eta(x) \leq (4d)^n] \leq  1-(\log 4d)^{-\delta} n^{-\delta}.
\]
As the number of particles in that ball
is of order~$n^{\theta d}$, one
gets that with probability at least
$1-\exp\{-C_2n^{\theta d-\delta}\}$ at time~$n^\theta$ one will
have at least one activated site~$x$ with $\eta(x)\geq (4d)^n$
and $\|x\|\leq C_1n^\theta$.
 Note the following simple fact:
If~$x$ contains at least $(2d)^n$ active particles and
$\|x-y\|_1\leq n$,
then until time~$n$ with probability bounded away from~$0$ at
least one of those particles will hit~$y$. Using this fact,
as~$x$ really contains at least~$2^n$ groups of~$(2d)^n$
particles, we get
that with overwhelming probability all the particles in the diamond
${\mathfrak D}_{n-n^\theta}$ will be awake at time $n^\theta+n$, which
completes the proof of Theorem~\ref{full_diamond}. \qed

\section{Remarks about continuous time}
A continuous-time version of the frog model can also
be considered. Here we would like to
remark that in the continuous-time context
and for the case of bounded~$\eta$,
Theorem~\ref{maintheo} also holds and its proof
can be obtained by following the steps of our
proof for the discrete case. The difficulty
that arises is that, for continuous time, it is not
evident that $\mu(x)$ (defined
by~(\ref{limtint})) is strictly
positive for $x\neq 0$, i.e., we must rule out the
possibility that the continuous-time frog model
grows faster than linearly. To overcome that difficulty,
note the following fact (compare with Lemma~9 on page~16 of
Chapter~1 of~\cite{D}): there exist a positive number~$\beta$ such
that, being $\|x\|_1\geq \beta n$, $\bbP[T(0,x)<n]$ is exponentially
small in~$n$. This fact in turn follows from a
domination of the frog model by branching random walk.
So, we conclude that for a.s.\ bounded~$\eta$ our
method works well in the continuous-time context too.

\section*{Acknowledgements}
This work was done when O.S.M.~Alves was visiting the Statistics
Department of the Institute of Mathematics and Statistics of the
University of S\~ao Paulo. He is thankful to the probability group
of IME-USP for hospitality.
O.S.M.~Alves thanks CAPES/PICD, F.P.~Machado and S.Yu.~Popov thank
CNPq (300226/97--7 and 300676/00--0), and K.~Ravishankar thanks
NSF (DMS 9803267) and FAPESP (01/01459--0) for financial support.

\end{document}